\documentclass[10pt]{article}
\usepackage{graphicx}
\usepackage{amsmath}
\usepackage[dvips]{epsfig}
\usepackage{amssymb}

\begin{document}

\title { Analytic Continuation of $q$-Euler numbers and polynomials}
\author{Taekyun Kim$^{1}$  \\[0.5cm]
$^{1}$ School of Electrical Engineering and Computer Science ,\\
        Kyungpook  National University, Taegu 702-701, S. Korea\\
          {\it e-mail: tkim$@$knu.ac.kr }\\ }

\date{}
\maketitle

 {\footnotesize {\bf Abstract}\hspace{1mm}
 In this
paper we study  that the $q$-Euler numbers and polynomials
 are analytically continued to $E_q(s)$. A new formula for the Euler's $q$-Zeta function $\zeta_{E,q}(s)$
 in terms of nested series of $\zeta_{E,q}(n)$
is derived. Finally we introduce the new concept of the dynamics of
analytically  continued $q$-Euler numbers and  polynomials. }



\section{Introduction }

 Throughout this paper, $\mathbb{Z}, \mathbb{R}$ and  $\mathbb{C}$ will denote the ring of  integers,
  the field of real numbers and the complex numbers, respectively.

When one talks of $q$-extension, $q$ is variously considered as an
indeterminate, a complex numbers or $p$-adic numbers. Throughout
this paper, we will assume that $q\in\mathbb{C} $ with $|q |<1$. The
$q$-symbol $[x]_q $ denotes  $[x]_q =  { {1- q^x} \over  {1-q }},$
(see [1-16]).

In this paper we study  that the $q$-Euler numbers and polynomials
 are analytically continued to $E_q(s)$. A new formula for the Euler's $q$-Zeta function $\zeta_{E,q}(s)$
 in terms of nested series of $\zeta_{E,q}(n)$
is derived. Finally we introduce the new concept of the dynamics of
analytically  continued $q$-Euler numbers and polynomials.

\section{  Generating $q$-Euler polynomials and numbers}

 For $h \in \mathbb{Z}$, the $q$-Euler polynomials were defined as
 $$ \sum_{n=0}^\infty {{E_n(x,h |q)} \over {n!}} t^n = [2]_q \sum_{n=0}^\infty (-1)^n q^{hn} e^{[n+x]_qt},\eqno (2.1) $$
 for $ x, q \in \mathbb{C}, $ cf. [1,7].
 In the special case $x=0$, $E_n(0,h |q)=E_n(h|q)$ are called
$q$-Euler numbers, cf. [1,2,3,4]. By (2.1), we easily see that
$$ E_n(x,h |q)={{[2]_q}\over{(1-q)^n}} \sum_{l=0}^n \binom nl (-1)^l {{1} \over {1+q^{l+h}}} q^{lx},
 \;\;  cf.   [7,8], \eqno (2.2) $$
 where $\binom nj $ is binomial coefficient.
From (2.1), we derive
 $$  E_{n,q}(x,h|q)=(q^x E(h|q)+[x]_q)^n
$$
with the usual convention of replacing $E^n(h|q)$ by $E_n(h|q).$
In the case $h=0$, $E_n(x,0|q)$ will be symbolically written as
$E_{n,q}(x).$
Let $G_q(x,t)$ be generating function of $q$-Euler
polynomials as  follows:
$$G_q(x,t)=\sum_{n=0}^{\infty} E_{n,q}(x)\frac{t^n}{n!}. \eqno(2.3)$$
Then we easily see that
$$G_q(x,t)= [2]_q  \sum_{k=0}^\infty (-1)^k e^{[k+x]_qt}. \eqno(2.4)$$
For $x=0,  E_{n,q}=E_{n,q}(0)$ will be called $q$-Euler numbers.

 From (2.3), (2.4),  we easily derive the following:
 For $k$(= even) and $n \in \mathbb{Z}_+$, we have
 $$ E_{n,q}{(k)}-E_{n,q}= [2]_q  \sum_{l=0}^{k-1} (-1)^l [l]_q^{n}. \eqno(2.5)$$
 For $k$(= odd) and $n \in \mathbb{Z}_+$, we have
$$ E_{n,q}{(k)}+E_{n,q}= [2]_q  \sum_{l=0}^{k-1} (-1)^l [l]_q^{n}.\eqno(2.6)$$
By (2.4), we easily see that
$$ E_{m,q}(x)= \sum_{l=0}^m \binom ml q^{xl} E_{l,q} [x]_q^{m-l}. \eqno(2.7)$$
From (2.5), (2.6), and (2.7), we derive
$$[2]_q \sum_{l=0}^{k-1} (-1)^{l-1} [l]_q^n = (q^{kn}-1) E_{n,q} + \sum_{l=0}^{k-1} \binom nl q^{kl} E_{l,q} [k]_q^{n-l}, \eqno(2.8)$$
where $k $(= even) $\in \mathbb{N}$. For $k$(= odd) and $n \in
\mathbb{Z}_+$, we have
$$[2]_q \sum_{l=0}^{k-1} (-1)^{l} [l]_q^n = (q^{kn}+1) E_{n,q} + \sum_{l=0}^{k-1} \binom nl q^{kl} E_{l,q} [k]_q^{n-l}. \eqno(2.9)$$

\section{  $q$-Euler zeta function }
It was known that the Euler polynomials are defined as
$$ \dfrac{2}{e^t+
1}e^{xt}= \sum_{n=0}^\infty  \dfrac{E_n(x)}{n!} t^n, \quad  |t| <
\pi, \text{ cf. [1-16].} \eqno(3.1)$$ For $s \in \mathbb{C}, x\in
\mathbb{R}$ with $0 \leq x<1$, define
$$ \zeta_E(s,x)= 2\sum_{n=0}^\infty {{(-1)^{n}}\over{(n+x)^{s}}}, \mbox{ and  }
\zeta_E(s)= 2\sum_{n=1}^\infty {{(-1)^{n}}\over{n^{s}}}. \eqno
(3.2)$$ By (3.1) and (3.2) we see that Euler numbers are related to
the Euler zeta function as
$$\zeta_E(-n)=E_n, \quad \zeta_E(-n,x)=E_n(x).$$
For $s, q, h \in \mathbb{C}$ with $|q|<1$, we define $q$-Euler zeta function as follows:
$$ \zeta_{E,q}(s,x|h)= [2]_q \sum_{n=0}^\infty {{(-1)^{n}q^{nh}}\over{[n+x]_q^{s}}},
\mbox{ and  } \zeta_{E,q}(s|h)= [2]_q \sum_{n=1}^\infty
{{(-1)^{n}q^{nh}}\over{[n]_q^{s}}}. \eqno (3.3)$$ For $k \in
\mathbb{N}, h \in \mathbb{Z}, $ we have
$$ \zeta_{E,q}(-n|h)= E_n(h |q). $$
In the special case $h=0$, $ \zeta_{E,q}(s|0)$ will be written
as $ \zeta_{E,q}(s) $. For $s \in \mathbb{C}$,  we note that
$$\zeta_{E,q}(s)=[2]_q \sum_{n=1}^\infty
{{(-1)^{n}}\over{[n]_q^{s}}}. $$

We now consider the function $E_q(s)$ as the analytic continuation of Euler numbers.
All the $q$-Euler numbers $E_{n,q}$ agree with $E_q(n)$, the analytic continuation of Euler numbers evaluated at $n$,
$$E_q(n)=E_{n,q} \text{ for } n \geq 0.$$
Ordinary Euler numbers are defined by
$$\dfrac{2}{e^t+1}= \sum_{n=0}^\infty E_n \dfrac{t^n}{n!}, \quad |t| < \pi. \eqno (3.4)$$
By (3.4),  it is easy to see that
$$E_0=1, \text{ and }
E_n=-\frac{1}{2}\sum_{l=0}^{n-1}\binom nl E_l, \text{
$n=0,1,2,\cdots .$}$$
 From (2.9) and (3.3), we can consider the $q$-extension
of Euler numbers $E_n$ as follows:
$$E_{0,q}=\frac{[2]_q}{2}, \text{ and }E_{n,q}=-\dfrac{1}{[2]_{q^{n}}} \sum_{l=0}^{n-1} \binom nl q^l
E_{l,q}, n=1,2,3, \cdots, \eqno(3.5)$$

In fact, we can express $E_q'(s)$ in terms of $\zeta_{E,q}'(s)$, the
derivative of $\zeta_{E,q}(s)$.
$$E_q(s)=\zeta_{E,q}(-s), E_q'(s)=\zeta_{E,q}'(-s),  E_q'(2n+1)=\zeta_{E,q}'(-2n-1), \eqno(3.6)$$
for $n \in \mathbb{N} \cup \{ 0 \}$.
This is just the differential of the functional equation and so verifies the consistency of
$E_q(s)$ and   $ E_q'(s)$ with $E_{n,q}$ and $\zeta(s)$.

From the above analytic continuation of $q$-Euler numbers, we derive
$$ \aligned
& E_q(s)=\zeta_{E,q}(-s),  E_q(-s)=\zeta_{E,q}(s) \\
& \Rightarrow E_{-n,q}=E_q(-n)=\zeta_{E,q} (n), n \in \mathbb{Z}_+.
 \endaligned \eqno(3.7) $$
The curve $E_q(s) $ runs through the points $E_{-n,q}$ and grows $
\sim n$ asymptotically as $(-n)\rightarrow -\infty$. The curve
$E_q(s)$ runs through the point $E_q(-n)$ and $ \lim_{n \to \infty}
E_q{(-n)}= \lim_{n \to \infty}\zeta_{E,q}(n)=-2 .$ From (3.5), (3.6)
and (3.7), we note that

$$\zeta_{E,q}(-n)=E_q(n)\mapsto \zeta_{E,q}(-s)=E_q(s).$$

\section{ Analytic continuation of  $q$-Euler polynomials }
For consistency with the redefinition of $ E_{n,q} = E_q(n)$ in
(4.5) and (4.6), we have
$$E_{n,q}(x)= \sum_{k=0}^n \binom nk E_{k,q} q^{kx} [x]_q^{n-k}. $$
Let $\Gamma (s)$ be the gamma function. Then the analytic
continuation can be obtained as

$$ \aligned
& n \mapsto s \in \mathbb{R} , x \mapsto w \in \mathbb{C}, \\
& E_{k,q} \mapsto E_q(k+s-[s])= \zeta_{E,q}(
-(k+(s-[s]))), \\
& \binom nk \mapsto
 {{\Gamma(1+s)}\over{\Gamma(1+k+(s-[s]))\Gamma(1+[s]-k)}}\\
& \Rightarrow E_{n,q}(s)\mapsto E_q(s,w )  = \sum_{k=-1}^{[s]}
{{\Gamma (1+s) E_q(k+s-[s])
q^{(k+s-[s])w}[w]_q^{[s]-k}}\over{\Gamma(1+k+(s-[s])) \Gamma
(1+[s]-k)}} \\
& = \sum_{k=0}^{[s]+1} {{\Gamma (1+s) E_q((k-1)+s-[s])
q^{((k-1)+s-[s])w}[w]_q^{[s]+1-k}}\over{\Gamma(k+(s-[s])) \Gamma
(2+[s]-k)}},
 \endaligned
$$
where $[s]$ gives the integer part of $s$, and so $s-[s]$ gives
the fractional part.

Deformation of the curve $E_q(2,w)$ into the curve of
$E_q(3,w)$ via the real analytic continuation $E_q(s,w), 2
\leq s \leq 3, -0.5 \leq w \leq 0.5.$

ACKNOWLEDGEMENTS. This paper is supported by Jangjeon Mathematical
Society and Jangjeon Research Institute for Mathematical Science and
Physics( 2007-001-JRIMS 1234567)

\end{document}